\documentclass[11pt, reqno]{amsart} 
\usepackage[letterpaper, margin = .8in]{geometry} 
\usepackage{lmodern} 
\usepackage{mathrsfs} 

\usepackage{amsmath, amsthm, amssymb, amsfonts} 
\usepackage{mathtools} 
\usepackage{booktabs} 
\usepackage{enumitem} 
\usepackage{comment}
\usepackage{graphicx} 
\usepackage{caption}
\usepackage{subcaption}
\usepackage{thmtools}
\usepackage{thm-restate}
\usepackage[dvipsnames]{xcolor} 
\usepackage{color}
\usepackage{hyperref}
\usepackage{url}
\usepackage{breakurl}
\usepackage{float}
\usepackage{bbm}
\newcommand{\bburl}[1]{\textcolor{blue}{\url{#1}}}

\newcommand{\R}{\mathbb{R}}

\newcommand{\Z}{\mathbb{Z}}

\newcommand{\GL}{\mathrm{GL}}
\newcommand{\supp}{\mathrm{supp}}
\newcommand{\sym}{\mathrm{sym}}
\newcommand{\1}{\mathbbm{1}}

\newtheorem{thm}{Theorem}[section]
\newtheorem{lem}[thm]{Lemma}
\newtheorem{conj}[thm]{Conjecture}
\newtheorem{prop}[thm]{Proposition}

\newtheorem{defi}[thm]{Definition}
\newtheorem{rek}[thm]{Remark}
\newtheorem{hyp}[thm]{Hypothesis}

\newcommand{\nc}{\newcommand}
\nc{\cF}{\mathcal F}
\nc{\DD}{\mathbb D}
\nc{\TT}{\mathbb T}
\nc{\EE}{\mathbb E}
\nc{\Symp}{\mathsf{Sp}}
\nc{\SpOrthO}{\mathsf{SO(odd)}}
\nc{\SpOrthE}{\mathsf{SO(even)}}
\nc{\Orth}{\mathsf O}
\nc{\Unit}{\mathsf U}
\nc{\UnitSp}{\mathsf{USp}}

\usepackage[backend=biber, doi=false, style=alphabetic, maxalphanames=10, maxnames=10]{biblatex} 

\numberwithin{equation}{section}

\title[Extending support for $1$ and $2$-level densities under cancellation hypotheses]{Extending the support of $1$- and $2$-level densities for cusp form $L$-functions under square-root cancellation hypotheses}

\author{Annika Mauro}
\address{Department of Mathematics, Stanford University}
\email{amauro@stanford.edu}

\author{Jack B. Miller}
\address{Department of Mathematics, Yale University}
\email{jack.miller.jbm82@yale.edu}

\author{Steven J. Miller}
\address{Department of Mathematics and Statistics, Williams College}
\email{sjm1@williams.edu}

\date{\today}

\thanks{This work was supported by NSF grant DMS1947438 and Williams College, and is gratefully dedicated to Henryk Iwaniec in celebration of his 75\textsuperscript{th} birthday for his years of direct mentoring and friendship to the third named author, and his indirect to the other authors (and numerous other students of the third author).}



\begin{document}

\maketitle

\begin{abstract}
The Katz-Sarnak philosophy predicts that the behavior of zeros near the central point in families of $L$-functions agrees with that of eigenvalues near 1 of random matrix ensembles. Under GRH, Iwaniec, Luo and Sarnak showed agreement in the one-level densities for cuspidal newforms with the support of the Fourier transform of the test function in $(-2, 2)$. They increased the support further under a square-root cancellation conjecture, showing that a $\GL(1)$ estimate led to additional agreement between number theory and random matrix theory. We formulate a two-dimensional analog and show it leads to improvements in the two-level density. Specifically, we show that a square-root cancellation of certain classical exponential sums over primes increases the support of the test functions such that the main terms in the $1$- and $2$-level densities of cuspidal newforms averaged over bounded weight $k$ (and fixed level $1$) converge to their random matrix theory predictions. We also conjecture a broad class of such exponential sums where we expect improvement in the case of arbitrary $n$-level densities, and note that the arguments in \cite{ILS} yield larger support than claimed.
\end{abstract}


\tableofcontents

\section{Introduction}

\subsection{Background}

Since the observations of Montgomery and Dyson 50 years ago, random matrix theory has provided a guide to predicting the behavior of quantities related to the zeros and values of $L$-functions; we focus on the behavior of zeros here. Initially this agreement was limited to theoretical results on the pair correlation of zeros of the Riemann zeta function, and then extended to include the $n$-level correlations of automorphic forms and numerical results on the spacings between zeros; see \cite{Ga, Hej, Mon, RS, Od1, Od2}. These statistics concern the behavior high up on the critical line, and are thus insensitive to finitely many zeros. In particular, they miss the effects from zeros at the central point, which are of great importance in a variety of number theory problems, from the Birch and Swinnerton-Dyer Conjecture and the ranks of the group of rational solutions of elliptic curves to bounding class numbers; see \cite{BSD1, BSD2, Go}.

To remedy this, Katz and Sarnak \cite{KS1, KS2} introduced new statistics, the $n$-level densities; there is now an extensive literature on agreement between number theory and random matrix theory here; see \cite{AAILMZ, AM, C--, CI, DHKMS1, DHKMS2, DM1, DM2, ER-GR, FiM, FI, Gao, GK, GJMMNPP, Gu, HM, HR1, HR2, ILS, LM, Mil1, MilPe, OS1, OS2, RR, Ro, Rub, ShTe, Ya, Yo1, Yo2} and the references therein, as well as \cite{BFMT-B, Con, FM} for more on the history of the connections between the two subjects.

This work is a continuation of the seminal paper \cite{ILS}; in brief, we extend their results for certain one-level densities to two-level, and discuss how to generalize to arbitrary $n$. Specifically, we show that a natural generalization of their Hypothesis S on cancellation in certain prime sums, which led to increasing the support for the one-level density, leads to increased support where the $2$-level density of certain families of cuspidal newforms and random matrix theory agree. Such results have applications in bounding the order of vanishing at the central point. We assume the reader is familiar with the basics of $L$-functions; see for example \cite{IK} for details. The following definitions are standard; we paraphrase from \cite{LiM} as we will use their framework to convert our results on increased support to estimates on weighted order of vanishing.

A function $\Phi:\R^n\to\R$ is Schwartz if it is infinitely differentiable and it and all of its derivatives decay faster than any polynomial. In our setting, an even Schwartz function $\Phi$ with compactly supported Fourier transform is called a \emph{test function}. We frequently assume the Generalized Riemann Hypothesis (GRH) holds for each $L(s, f)$, and write the non-trivial zeros of a cuspidal newform $L(s, f)$ of level $N$ and weight $k$ by \[ \rho^{(j)}_f \ = \  \frac12 + i \gamma_f^{(j)} \] for $\gamma_f^{(j)} \in \R$ increasingly ordered and centered about zero.\footnote{In many of the calculations of Bessel-Kloosterman terms, we need GRH for Dirichlet $L$-functions. We can often avoid GRH for cuspidal newforms at the expense of more involved calculations, but if GRH fails while we can still formally calculate these statistics, as the zeros are no longer on the critical line we lose the correspondence with physics.} The number of zeros with $|\gamma_f^{(j)}|$ bounded by an absolute large constant is of order $\log c_f$ for some constant $c_f > 1$; this is known as the \emph{analytic conductor}.

\begin{defi}[$n$-Level Density]
The $n$-level density of an $L$-function $L(s,f)$ with respect to a test function $\Phi:\R^n\to\R$ with compactly supported Fourier transform is defined as
\begin{equation}
    D_n (f; \Phi) \ := \  \sum_{\substack{j_1, \dots, j_n \\ j_i \neq \pm j_k}} \Phi \left( \frac{\log c_f}{2\pi} \gamma_f^{(j_1)}, \dots, \frac{\log c_f}{2\pi} \gamma_f^{(j_n)} \right).	\label{def:density}
\end{equation}
Note that the existence of the $n$-level density does not depend on GRH.
\end{defi}

One of the most important applications of the $n$-level density is to obtaining bounds on the order of vanishing to a given order at the central point. For such results, we need $\Phi$ to be non-negative and positive at the origin.

Unlike the $n$-level correlations, the sum \eqref{def:density} is hard to study for an individual $f$ because by choice of $\Phi$ it essentially captures only a bounded number of zeros. Thus we study averages over finite subfamilies $\cF(Q) := \{ f \in \cF : c_f \in I(Q)\}$ (which are parametrized by some quantity $Q$ such that as $Q$ tends to infinity, the size of the subfamily tends to infinity as well), namely
\begin{equation}
    \EE[D_n (f; \Phi), Q] \ := \  \frac{1}{\# \cF (Q)} \sum_{f \in \cF (Q)} D_n(f; \Phi). \label{eq:averagessecond}
\end{equation} Common choices are $I(Q) = \{Q\}, \{1, 2, \dots, Q\}$ and $\{Q, Q+1, \dots, 2Q\}$.

The Katz-Sarnak density conjecture \cite{KS1,KS2} asserts that the $n$-level density of a family of $L$-functions, in the limit as the conductors tend to infinity, converges to the $n$-level density of eigenvalues of a classical compact group as the matrix sizes tend to infinity. Explicitly,  if $\cF$ is a ``good'' family of $L$-functions and $\Phi$ is not zero at the origin, then there exists a distribution $W_{n, \cF}$ such that
\begin{eqnarray}
    & & \lim_{Q \to \infty} \EE[D_n (f; \Phi), Q] \ = \  \frac{1}{\Phi(0, \dots, 0)} \int_{\R^n} \Phi(x)\cdot W_{n, \cF} (x) \ dx_1 \cdots dx_n.
\end{eqnarray}

We have the following expansions for the quantities above; though for computational purposes it is often advantageous to use an alternative expansion due to Hughes-Miller \cite{HM}, which writes the $n$-level density as a sum of terms emerging as the support increases.

\begin{thm}[Determinant Expansion \cite{KS1}]\label{thm:KSdet} Let $K(x) = \frac{\sin{\pi x}}{\pi x}$ and $K_{\epsilon}(x,y) = K(x-y)+\epsilon K(x+y)$. Then the $n$-level densities have the following distinct closed form determinant expansions for each corresponding symmetry group:
\begin{align}
	W_{n, \SpOrthE} (x)
		&\ = \  \det \left( K_1 (x_i, x_j) \right)_{i, j \leq n}, \label{eq:nlevelSOeven} \\
	W_{n, \SpOrthO} (x)
		&\ = \  \det \left( K_{-1} (x_i, x_j) \right)_{i, j, \leq n} + \sum_{k  =  1}^n \delta_0 (x_k) \det \left( K_{-1} (x_i, x_j) \right)_{i, j, \neq k},  \label{eq:nlevelSOodd}\\
	W_{n, \Orth} (x)
		&\ = \  \frac12 W_{n, \SpOrthE} (x) + \frac12 W_{n, \SpOrthO} (x), \label{eq:nlevelSpOrthE}\\
	W_{n, \Unit} (x)
		&\ = \  \det \left( K_0 (x_i, x_j) \right)_{i, j, \leq n}, \label{eq:nlevelUni}\\
	W_{n, \Symp} (x)
		&\ = \ \det \left( K_{-1} (x_i, x_j) \right)_{i, j, \leq n}.\label{eq:nlevelSymp}
\end{align}
\end{thm}

As remarked, one of the main applications of this statistic is to bound the order of vanishing of a family of $L$-functions at the central point by choosing a test function which is non-negative and positive at the origin, see \cite{BCDMZ, DM, Fr, FrMil, LiM}. Many of these papers are concerned with trying to find the optimal test function for a given support, but already in the work of \cite{ILS} one sees that there may be only a negligible improvement in bounds from using the optimal functions derived from Fredholm theory over simple test functions. Thus, while there have been some recent advances in determining the optimal function for a given support, it has been more productive to increase the support and the level $n$ studied; however, as $n$ increases while the bounds obtained are better for the percentage of forms vanishing to order at least $r$ when $r$ is large, they are worse for small $r$. Thus there is a balancing act, with most of the effort devoted to finding the largest support possible, and then given that determining the best bounds by using easy to compute test functions for each $n$. We concentrate on increasing the support for certain families of cusp forms, which can then be fed into the machinery from Dutta-Miller \cite{DuM} to yield improved estimates.


In \cite{ILS}, the authors introduce the following hypothesis, whose implications are striking, allowing them to break $(-2, 2)$ for the support for certain families.

\begin{hyp}[Hypothesis S]
\label{hyp:HypothesisS}
Let $e(z) := e^{2\pi iz}$, $c$ be a positive integer and $a$ be an arbitrary residue class mod $c$. Let $H_1(\alpha,A)$ denote a $2$-parameter family of hypotheses, where $\alpha\in[1/2, 3/4]$, $A\in[0,\infty)$, each of which states that
\begin{equation}
\label{eq:HypothesisS}
H_1(\alpha,A): \qquad
\underset{p \equiv a (c)}{\sum_{p \le x}}
e\left(\frac{2\sqrt{p}}{c}\right)
\ \ll_\varepsilon \
c^A x^{\alpha+\varepsilon}
\end{equation}
holds uniformly over $c>0$ and residue classes $a$ mod $c$.
\end{hyp}

\subsection{Previous Work}

The reason why \cite{ILS} is interested in this family of hypotheses is because they increase the support of test functions for which the one-level density agrees with random matrix theory. We first review their result, and then discuss our generalization.

For convenience, as is often done in the subject, we assume our test function $\Phi$ is the product of one-dimensional test functions $\phi_i$. Below we confine our study to cusp forms of level 1 and weight $k$. Following \cite{ILS}, one is able to get better results on the support by averaging over $k$. This allows us to exploit some oscillation in the Bessel function factors that emerge in application of the Petersson formula. Note that rather than sum over all forms equally, each form $f$ of weight $k$ is weighted by the slowly varying factor $1/L(1, \sym^2(f))$ \cite{HL,I2}. There is a trade-off in studying this modification of the $n$-level densities; these harmonic weights arise naturally in the Petersson formula, and their introduction simplifies calculations below. Unfortunately their presence means that we cannot obtain results on bounding the order of vanishing at the central point, but instead obtain results on weighted vanishing. In many families these weights can be removed through additional work; see \cite{ILS}.

\begin{defi}[Agrees with Orthogonal RMT]
We say that a test function $\Phi$
\emph{agrees with Orthogonal Random Matrix Theory} if
\begin{equation}
\lim_{K\to\infty} \frac{1}{K} \ \sum_{k \le K} \ \frac{4\pi^2}{k-1}
\sum_{f\in H_k} \frac{1}{L(1,\sym^2(f))} D_n(f;\Phi)
\ \ = \  \
\int_{\R^n} \Phi(x)
\cdot W_{n,\Orth}(x)
\ dx_1\cdots dx_n,
\end{equation}
where $W_{n,\Orth}(x)$ is given by equation (\ref{eq:nlevelSpOrthE}), and $H_k = H_k^{\star}(1)$ are Hecke eigencuspforms of weight $k$ and level $1$, with normalization so that $a_f(1)=1$. The $L$-function $L(s,f)$ has coefficients $\lambda_f(n) \coloneqq a_f(n) n^{-(k-1)/2}$ so that $\lambda_f(p)\in[-2,2]$ and the $L$-function $L(s,f)$ is symmetric about $\mathfrak{Re}(s) = 1/2$.
\end{defi}

\begin{rek} For the family considered above, previous results show that the underlying symmetry group is orthogonal, hence our comparison with the orthogonal behavior. For other families we would just use the corresponding densities from Theorem \ref{thm:KSdet}. \end{rek}

We are now able to state the following result from \cite{ILS}, which extends the support for the family of level 1 cusp forms to beyond $(-2, 2)$.

\begin{thm}[One-level extended support]
\label{thm:HypothesisS}
Assume GRH. Then $H_1(\alpha,A)$ implies that a test function $\phi$ agrees with random matrix theory if $\supp(\widehat{\phi}) \subset (-\sigma,\sigma)$ for $\sigma = \min\{5/2, 2 + (6-8\alpha)/(1+2A+4\alpha)\}$. More specifically, we show for all $h\in C_0^\infty(\R^+)$ with $\widehat{h}(0) \neq 0$ that
\begin{equation}
\label{eq:1level Theorem Averaging}
\lim_{K\to \infty}\frac{1}{\widehat{h}(0) K} \
\sum_{k \ {\rm even}} \
\frac{4\pi^2}{k-1} \
h\left(\frac{k-1}{K}\right)
\sum_{f \in H_k}
\frac{1}{L(1,\sym^2(f))} D_1(f;\phi) \ \to \
\int_{\R} \phi(x) \cdot W_{1,\Orth}(x) \,dx.
\end{equation}
\end{thm}

\subsection{New Results}

In their manuscript from 2000, Iwaniec, Luo and Sarnak \cite{ILS} state and prove the Theorem \ref{thm:HypothesisS} for $\sigma = 2+(12-16\alpha)/(5+4A+8\alpha)$, which is less than $2 + (6-8\alpha)(1+2A+4\alpha)$; however, a careful analysis of their arguments give more than what they claimed. In Section \ref{sec:HypothesisS}, we give a proof of Theorem \ref{thm:HypothesisS} based on Section 10 of \cite{ILS}, and discuss in Appendix \ref{sec:HistoryOfThm} the difference in the claimed support and what the calculation yields.

Our main theorem, Theorem \ref{thm:HypothesisT}, is a natural extension of Theorem \ref{thm:HypothesisS}, showing that square-root cancellation hypotheses also extend to the case of the $2$-level density. For ease of exposition and to highlight the issues, we focus on the $n=2$ case, though similar calculations should hold in general.

\begin{hyp}\label{hyp:HypothesisT}
Let $e(z) := e^{2\pi iz}$, $c$ be a positive integer, and $a_1,a_2$ be arbitrary residue classes mod $c$. Let $H_2(\alpha,A)$ denote a $2$-parameter family of hypotheses, where $\alpha\in[1/2, 3/4]$, $A\in[0,\infty)$, each of which states that
\begin{equation}
\label{eq:HypothesisT}
H_2(\alpha,A): \qquad
\underset{p_1 \equiv a_1 (c)}{\sum_{p_1 \leq x_1}} \ \
\underset{p_2 \equiv a_2 (c)}{\sum_{p_2 \leq x_2}}
e\left(\frac{2\sqrt{p_1p_2}}{c}\right)
\ \ll_\varepsilon \
c^A (x_1x_2)^{\alpha+\varepsilon}
\end{equation}
holds uniformly over $c>0$ and residue classes $a_1,a_2$ mod $c$.
\end{hyp}

It should be noted that hypothesis $H_2(\alpha,A)$ implies $H_1(\alpha,A)$, because one may take $x_1$ or $x_2$ constant to obtain a single sum over primes.
A natural extension of Theorem \ref{thm:HypothesisS} is our main theorem, proven in Section \ref{sec:HypothesisT}.

\begin{thm}[Two-level extended support]\label{thm:HypothesisT}
Assume GRH. Then $H_2(\alpha,A)$ implies that a test function $\Phi(x) = \phi_1(x_1)\phi_2(x_2)$ agrees with random matrix theory if $\supp(\widehat{\phi_i}) \subset (-\sigma_i,\sigma_i)$ for \\ $\sigma_1+\sigma_2 = 2 + (6-8\alpha)/(3+2A+4\alpha)$. More specifically, we show for all $h\in C_0^\infty(\R^+)$ with $\widehat{h}(0) \neq 0$ that
\begin{equation}\label{eq:HypothesisTtheoremeq}
\lim_{K\to\infty}\frac{1}{\widehat{h}(0) K} \
\sum_{k \ {\rm even}} \
\frac{4\pi^2}{k-1} \
h\left(\frac{k-1}{K}\right)
\sum_{f \in H_k}
\frac{1}{L(1,\sym^2(f))} D_2(f;\Phi) \ \to \
\int_{\R^2} \Phi(x) \cdot W_{2,\Orth}(x) \,dx_1dx_2.
\end{equation}
\end{thm}

Using the machinery from \cite{DuM} for an orthogonal ensemble with the naive test function we obtain the following immediately from Theorem \ref{thm:HypothesisT}.

\begin{thm}\label{thm:rankbounds} Taking $\phi$ to be the naive test function
\begin{equation}
\phi_{\rm Naive}(x) \ = \ \left(\frac{\sin(\pi v x)}{\pi vx}\right)^2,
\end{equation}
the \emph{weighted percentage of forms vanishing} to order at least $r$, denoted $P_r(\mathcal{F})$, is bounded by
\begin{equation}
P_r(\mathcal{F}) \ \leq \ \min\left\{
\frac{1}{r}\left(\frac{1}{2} + \frac{1}{v_1}\right),
\frac{1}{3(r-1/v_2-1/2)^2}
\right\}.
\end{equation}
Here $v_1$ and $v_2$ denote the available support for the $1$-level and $2$-level density respectively, and the form $f$ is weighted by the factor from \eqref{eq:HypothesisTtheoremeq} (normalized so that the sum of the weights equals 1). Assuming the strongest version of Hypothesis \ref{hyp:HypothesisT}, we may take $v_1 = 2.5$ and $v_2 = 1.2$. These values yield
\begin{equation}
P_1(\mathcal{F}) \ \leq \ 0.9000, \quad
P_2(\mathcal{F}) \ \leq \ 0.4500, \qquad
P_3(\mathcal{F}) \ \leq \ 0.1200, \qquad
P_4(\mathcal{F}) \ \leq \ 0.0469, \qquad
P_5(\mathcal{F}) \ \leq \ 0.0248.
\end{equation}
\end{thm}

\begin{rek} We quickly comment on the conditional nature of these results. In \cite{ILS} there are two places GRH is used. The first is for convenience to bound certain prime sums of cuspidal newforms, which they remark can be bypassed by additional appeals to the Petersson formula (see the comments after their equation (4.24)).\footnote{While we need GRH to have the zeros lie on the critical line, and thus have a direct comparison to eigenvalues of Hermitian matrices or energy levels of heavy nuclei, such an assumption is only needed for this interpretation or correspondence; the $n$-level densities exist whether or not GRH holds.} The second is when one splits by sign for square-free level $N$ tending to infinity, GRH is needed to analyze the main term contribution and size of the error of the Bessel-Kloosterman term; however, we only need this for Dirichlet $L$-functions. If we do not split by sign and consider level 1 cuspidal newforms for $k$ up to $K$ (appropriately weighted), then Hypothesis S, a ``${\rm GL}_1$'' exponential sum, suffices. \end{rek}

\begin{rek} We expect in the limit that 50\% of the forms in our family should be rank 0 and 50\% rank 1; thus (paraphrasing comments one of us heard from Iwaniec in graduate classes) for $r \ge 2$ we may interpret Theorem \ref{thm:rankbounds} as providing better upper bounds on 0. \end{rek}

We first review the proof for the one-level density in \S\ref{sec:HypothesisS} to fix notation and then extend to the two-level in the next section. For notational convenience we prove Theorem \ref{thm:HypothesisT} in the case when $\phi_1=\phi_2$; a similar argument holds in general (see for example the arguments in \cite{LiM}, which extend the $n$-level density results of \cite{HM} from identical test functions to the more general case). We then discuss further generalizations to larger $n$.

\section{The One-Level Case}\label{sec:HypothesisS}

We follow the arguments in \cite{ILS} to prove Theorem \ref{thm:HypothesisS}, but with extra care towards the step that derives their equation (10.17), since the expression in their original manuscript differs slightly from our estimation. This also sets the notation we need for our two-level result, as well as isolating certain one-level results that are used again.

We begin by justifying normalizing the one-level density sum
\begin{equation}
\label{eq:1levelsum}
\mathscr{B}(K) \ \coloneqq \
\sum_{k \ {\rm even}} \
\frac{4\pi^2}{k-1} \
h\left(\frac{k-1}{K}\right)
\sum_{f \in H_k}
\frac{1}{L(1,\sym^2(f))} D_1(f;\phi)
\end{equation}
by a factor of $\widehat{h}(0)K$. That is, given the total weighting
\begin{equation}
B(K) \ \coloneqq \
\sum_{k \ {\rm even}} \
\frac{4\pi^2}{k-1} \
h\left(\frac{k-1}{K}\right)
\sum_{f \in H_k}
\frac{1}{L(1,\sym^2(f))},
\end{equation}
we show that $B(K) = \widehat{h}(0)K + O(1)$, and thus either gives the same main terms as $K\to\infty$. This follows from a special case of the Petersson trace formula:
\begin{equation}
B(K) \ = \
\sum_{k \ {\rm even}} \
2h\left(\frac{k-1}{K}\right) \
\Delta_k(1,1), \qquad
\Delta_k(1,1) \ = \
\frac{2\pi^2}{k-1} \sum_{f\in H_k} \frac{1}{L(1,\sym^2(f))},
\end{equation}
where the so-called \emph{trace} $\Delta_k(m,n)$ is defined to be
\begin{equation}
\label{eq:Petersson Trace}
\Delta_k(m,n) \coloneqq
\sum_{f\in H_k} \frac{\Gamma(k-1)}{(4\pi\sqrt{mn})^{k-1}} \frac{a_f(m) \ a_f(n)}{\langle f,f \rangle}
\ = \
\frac{2\pi^2}{k-1} \sum_{f\in H_k} \frac{\lambda_f(m)\lambda_f(n)}{L(1,\sym^2(f))}
\end{equation}
where $\langle\cdot,\cdot\rangle$ is the Petersson inner product and $a_f$ are the Fourier coefficients of the cusp form $f$.

\begin{prop}[Petersson trace formula] Let $\delta(\cdot,\cdot)$ be Kronecker's delta, $J_{k-1}$ a Bessel function of the first kind, and
\begin{equation}
S(m,n;c) \ \coloneqq \
\sum_{d \ ({\rm mod}\ c)} \hspace{-12pt}^\star \ \
e\left(\frac{md+n\overline{d}}{c}\right)
\end{equation}
is a Kloosterman sum, where $\sum\hspace{0pt}^\star$ denotes summing over primitive residue classes. Then the Petersson formula is
\begin{equation}
\label{eq:Petersson formula}
\Delta_k(m,n) \ = \
\delta(m,n) +
2\pi i^k \
\sum_{c = 1}^\infty \
\frac{S(m,n;c)}{c}
J_{k-1}\left(\frac{4\pi\sqrt{mn}}{c}\right).
\end{equation}
\end{prop}

Careful estimation of the Petersson trace formula for special values as in Corollary 2.3 of \cite{ILS} gives $\Delta_k(1,1) = 1 + O(2^{-k})$, and so our claim about $B(K)$ follows by interpreting the main term $H = \sum_{k \ {\rm even}} 2h\left(\frac{k-1}{K}\right)$ as a Riemann sum of the function $2h(x)$ divided by the mean spacing $2/K$, hence $H = \widehat{h}(0)K + O(1)$ where the implied constant depends on $h$, e.g.\ using the fact that $h$ has bounded derivative (though the weaker condition of $h$ having bounded variation would suffice).

In performing asymptotic analysis, we are interested in the weighted sum of traces
\begin{equation}
\label{eq:mathscrB(,)}
\mathscr{B}(m,n) \ \coloneqq \
\sum_{k \ {\rm even}} 2h\left(\frac{k-1}{K}\right) \Delta_k(m,n).
\end{equation}
Inserting (\ref{eq:1levelexplicit}) into (\ref{eq:1levelsum}) while taking $R \asymp K^2$ to be on order of the average conductor, we have
\begin{equation}
\mathscr{B}(K) \ = \
\widehat{h}(0) K\big<\phi,W_{1,\Orth}\big> - \mathscr{P}(\phi) + O\left(K\frac{\log\log K}{\log K}\right),
\end{equation}
where the implied constant depends on $\phi$, and $\mathscr{P}(\phi)$ is the weighted sum over local factors
\begin{equation}
\label{eq:mathscrP}
\mathscr{P}(\phi) \ \coloneqq \
\sum_{k \ {\rm even}} \
\frac{4\pi^2}{k-1} \
h\left(\frac{k-1}{K}\right)
\sum_{f \in H_k}
\frac{1}{L(1,\sym^2(f))}
\textbf{P}(f;\phi)
\ = \
\sum_p \ \mathscr{B}(p,1) \ \widehat{\phi}\left(\frac{\log p}{2\log K}\right) \frac{\log p}{p^{1/2}\log K}.
\end{equation}

It is at this stage the goal of the analysis becomes clear: estimate $\mathscr{P}(\phi)$ by extracting a main term and bounding the error. In the case of the one-level density, the prime sum $\mathscr{P}(\phi)$ does \emph{not} contribute to the main term when $\supp(\widehat{\phi}) \subset (-\sigma,\sigma)$ for some $\sigma$ we would like to determine, giving us agreement with random matrix theory. For example, without the use of Hypothesis \ref{hyp:HypothesisS}, \cite{ILS} immediately show using the Petersson trace formula and Weil's estimate that $\mathscr{B}(p,1) = O(p^{1/2}K^{-4})$, and so taking $\sigma = 2$ means that $p$ runs up to $P \ll K^{4-\delta}$ assuming $\supp(\widehat{\phi}) \subseteq [-2+\delta,2-\delta]$ for some positive $\delta$, which gives $\mathscr{P}(\phi) \ll K^{-\delta}$.

To increase $\sigma$, it is necessary to deal with a worse error term that cannot be absorbed into $O(p^{1/2}K^{-4})$ when $\sigma > 2$. Namely, by applying the Petersson formula to $\mathscr{B}(m,n)$ and approximating the Bessel function sums using standard techniques (see Corollary 8.2 in \cite{ILS} and then the top of page 86 in \cite{I1}), one may derive the following.

\begin{lem}[Lemma 10.1 in \cite{ILS}]
\label{lem:ILSlemma10.1}
We have
\begin{align}
\label{eq:ILSlemma10.1}
\mathscr{B}(m,n) \ & = \
\widehat{h}(0)K\delta(m,n) + O(\delta(m,n) + \sqrt{mn} K^{-4}) \\ & \qquad
- \pi^{1/2}(mn)^{-1/4} K \ {\rm Im}\left(
\overline{\zeta}_8 \ \sum_{c=1}^\infty \
c^{-1/2} S(m,n;c) \  e\left(\frac{2\sqrt{mn}}{c}\right) \hslash\left(\frac{cK^2}{8\pi\sqrt{mn}}\right)
\right), \nonumber
\end{align}
where the implied constant depends on $h$, and the function $\hslash$ is defined to be the transform
\begin{equation}
\label{eq:hslash}
\hslash(v) \ \coloneqq \
\int_0^\infty
\frac{h(\sqrt{u})}{\sqrt{2\pi u}} e^{iuv}\,du.
\end{equation}
\end{lem}

The above lemma follows by applying the Petersson trace formula to (\ref{eq:mathscrB(,)}), and then estimating the weighted sum of Bessel functions $I(x) = \sum_{k \ {\rm even}} 2h\left(\frac{k-1}{K}\right) i^k J_{k-1}(x)$ (cf. the twisted character sum on page 86 of \cite{I1}). Crucial in the analysis is executing the sum over the weights $k$ to exploit the oscillation in the Bessel terms.

Using this expression for $\mathscr{B}(p,1)$ to estimate (\ref{eq:mathscrP}), we obtain by a simple triangle inequality
\begin{align}
\mathscr{P}(\phi) \ &\ll_h \
PK^{-4} \\ &\quad +
\left|\sum_p
K
\ {\rm Im}\left(
\overline{\zeta}_8
\sum_{c=1}^\infty
c^{-1/2} S(p,1;c) \ e\left(\frac{2\sqrt{p}}{c}\right) \hslash\left(\frac{cK^2}{8\pi\sqrt{p}}\right)
\widehat{\phi}\left(\frac{\log p}{2\log K}\right) \frac{\log p}{p^{3/4}\log K}\right)\right|, \nonumber
\end{align}
where $P$ is the largest prime such that $\frac{\log P}{2\log K}$ is greater than the support of $\widehat{\phi}$. Notice that we use a crude estimate such as $\sum_{p \leq P} 1 \ll P$ because logarithmic factors will not increase the support (as our results are for open and not closed intervals). We moved $\widehat{\phi}\left(\frac{\log p}{2\log K}\right)$ and $\frac{\log p}{p^{3/4}\log K}$ into the imaginary component, since they are both real.

Our next step is to interchange summing over $p$ versus $c$, and we replace $|{\rm Im}(\cdot)|$ by the full complex modulus $|\cdot|$ (we do not expect the real part to be significantly larger than the imaginary part, so this should not lead to a decrease in support).
\begin{align}
\label{eq:P(phi) bound before IBP}
\mathscr{P}(\phi) &\ll PK^{-4} +
K \ \sum_{c=1}^\infty
\ c^{-1/2}
\ \left| \sum_p \
S(p,1;c) e\left(\frac{2\sqrt{p}}{c}\right) \hslash\left(\frac{cK^2}{8\pi\sqrt{p}}\right)
\widehat{\phi}\left(\frac{\log p}{2\log K}\right) \frac{\log p}{p^{3/4}\log K}
\right|
\nonumber \\ &\ll
PK^{-4} +
K \ \sum_{c=1}^\infty
\ c^{-1/2} \
\sum_{a \ ({\rm mod } \ c)}\hspace{-12pt}^\star
\hspace{10pt}
\left| \hspace{2pt} \sum_{p \equiv a (c)} \
S(p,1;c) e\left(\frac{2\sqrt{p}}{c}\right) \hslash\left(\frac{cK^2}{8\pi\sqrt{p}}\right)
\widehat{\phi}\left(\frac{\log p}{2\log K}\right) \frac{\log p}{p^{3/4}\log K}
\right|
\nonumber \\ &=
PK^{-4} +
K \ \sum_{c=1}^\infty
\ c^{-1/2} \
\sum_{a \ ({\rm mod } \ c)}\hspace{-12pt}^\star
\hspace{10pt}
|S(a,1;c)|
\left| \hspace{2pt} \sum_{p \equiv a (c)} \
e\left(\frac{2\sqrt{p}}{c}\right) \hslash\left(\frac{cK^2}{8\pi\sqrt{p}}\right)
\widehat{\phi}\left(\frac{\log p}{2\log K}\right) \frac{\log p}{p^{3/4}\log K}
\right|.
\end{align}
We only sum over primitive residue classes mod $c$, each of which contains infinitely many primes, whereas the non-primitive classes consist of primes $p$ dividing $c$. Because $\hslash$ is rapidly decaying, only large primes $p\gg c^2K^4$ contribute, and so the non-primitive residue classes are absorbed into the Vinogradov notation.

We now perform summation by parts. For each sum over primes in a residue class, we put $\psi_c(x) = \hslash\left(\frac{cK^2}{8\pi\sqrt{x}}\right)\widehat{\phi}\left(\frac{\log x}{2\log K}\right)\frac{\log x}{x^{3/4}\log K}$ which is smooth and supported for primes $p < P$. By Abel summation, we obtain
\begin{equation}
\label{eq:AbelSummation1}
\sum_{2 \ \leq \ n \ \leq \ P} \
e\left(\frac{2\sqrt{n}}{c}\right)
\ \1_{n\in\{p \equiv a(c)\}} \
\psi_c(n)
\ = \
-\int_2^P \mathcal{E}_1(x) \psi_c'(x)\,dx,
\end{equation}
where
\begin{equation}
\mathcal{E}_1(x) \ \coloneqq \
\underset{p \le x}{\sum_{p \equiv a(c)}} \ e\left(\frac{2\sqrt{p}}{c}\right).
\end{equation}
Thus we see a natural opportunity to use Hypothesis \ref{hyp:HypothesisS} in our estimate for the local factors arising from the explicit formula applied to the weighted average of the one-level density $D_1(f;\phi)$.


From now on let us assume hypothesis $H_1(\alpha,A)$ found in (\ref{eq:HypothesisS}). We estimate $\psi_c'(x)$ as
\begin{equation}
\label{eq:psiprime}
\psi_c'(x) \ = \
\hslash'\left(\frac{cK^2}{8\pi\sqrt{x}}\right) \cdot O\left(\frac{cK^2}{x^{9/4}}\right) +
\hslash\left(\frac{cK^2}{8\pi\sqrt{x}}\right) \cdot O\left(\frac{1}{x^{7/4}}\right),
\end{equation}
where the implied constant depends on $\phi$.
We treat $\psi_c'(x)$ as $O(1)$ for small $c$ and as rapidly decaying for large $c$. The exact transition region for $c$ is governed by the argument of the rapidly decaying function $\hslash$, i.e., $cK^2/\sqrt{P} = \theta(1)$. Thus, we truncate our sum over $c$ in (\ref{eq:P(phi) bound before IBP}) at the value $C = P^{1/2} K^{\varepsilon-2}$.
Indeed, applying Weil's estimate and equation (\ref{eq:psiprime}) with $\hslash'(x),\hslash(x) = O_\Omega(x^{-\Omega})$ for $\Omega > 0$ large, we bound the tail of (\ref{eq:P(phi) bound before IBP}) by
\begin{align}
K\sum_{c \ \geq \ C}
c \int_{2}^P \mathcal{E}_1(x) \psi_c'(x)\,dx \ &\ll_{h,\Omega,\varepsilon} \
K\sum_{c \ \geq \ C}
c^{A+1} \int_{2}^P x^{\alpha+\varepsilon} \left(\frac{cK^2}{\sqrt{x}}\right)^{-\Omega} \cdot \left(\frac{cK^2}{x^{9/4}} + \frac{1}{x^{7/4}}\right)
\,dx \\ &\ll \
P^{\alpha+\varepsilon+\Omega/2} K^{1-2\Omega} C^{A+3-\Omega} \nonumber \\ &= \
K^{O(1)-\varepsilon\Omega} \nonumber.
\end{align}
Thus we have a balancing act, where we let $\varepsilon\to0$ while $O(1) - \varepsilon\Omega < 0$ to maximize the extension of support obtained by this method. This shows that we may neglect these terms for further computation, i.e., they are absorbed into the Vinogradov notation.

Focusing now on the small values of $c$, we estimate $\psi_c'(x) = O_{h}(K^{\varepsilon}x^{-7/4})$ for $x \geq c^2K^{4-2\varepsilon}$, handling the smaller values of $x$ using rapid decay estimates similar to above. This gives
\begin{align}
K \sum_{c \ \leq \ C} c \int_{c^2K^{4-2\varepsilon}}^P \mathcal{E}_1(x) \psi_c'(x)\,dx
\ &\ll_{h,\varepsilon} \
K^{1+\varepsilon} \sum_{c \ \leq \ C} c^{A+1} \int_{c^2K^{4-2\varepsilon}}^P x^{\alpha-7/4+\varepsilon} \,dx
\nonumber \\ &\ll \
P^{\alpha+A/2+1/4} K^{-2A-3+O(\varepsilon)}.
\end{align}
We usually write $O(\varepsilon)$ as $\varepsilon$, as this can be done by changing our definition of $\varepsilon$. Thus, we have derived
\begin{equation}
\label{eq:Final Estimate for P(phi)}
\mathscr{P}(\phi) \ \ll_{\phi,h,\varepsilon} \
PK^{-4} + P^{\alpha+A/2+1/4} K^{-2A-3+\varepsilon}.
\end{equation}
To understand an account of the various estimates for (\ref{eq:Final Estimate for P(phi)}), see Appendix \ref{sec:HistoryOfThm}.


Since $P = K^{2\sigma'}$ where $\sigma' = \frac{\log P}{2\log K}$ tends to the maximum of $\supp(\widehat{\phi})$, and so $\sigma'<\sigma$ for all large $K$. Since we have agreement with random matrix theory only if $\mathscr{P}(\phi) = o(K)$, each term in (\ref{eq:Final Estimate for P(phi)}) gives us an upper bound on our method for obtaining $\sigma$ in Theorem \ref{thm:HypothesisS}. Namely, $PK^{-4} = o(K)$ gives $\sigma \leq 5/2$, and $P^{\alpha+A/2+1/4}K^{-2A-3+\varepsilon}$ as $\varepsilon\to0$ gives $\sigma \leq 2 + (6-8\alpha)/(1+2A+4\alpha)$. This concludes our proof of Theorem \ref{thm:HypothesisS}.

The best support we could expect would be using hypothesis $H_1(\frac{1}{2},0)$, which gives an increased support of $\supp(\widehat{\phi}) \subset (-5/2,5/2)$; see Appendix \ref{sec:HistoryOfThm} for a history of this derivation.
Notice that $\sigma = 5/2$ can also be achieved using $H_1$ with respect to any $(\alpha,A)$ along the linear interpolation of $(0.5,0.5)$ and $(0.55,0)$.


\section{The Two-Level Case}\label{sec:HypothesisT}

The starting point is the explicit formula, and inclusion-exclusion to express the two-level density in terms of prime sums. The calculation is standard (see Section \ref{sec:Explicit Formula}) and yields equation (\ref{eq:2level explicit expansion}), which we restate below (recall $\textbf{P}(f;\phi)$ is defined in \eqref{eq:mathscrP}):
\begin{align}
D_2(f;\phi) \ &= \
\widehat{\phi}(0)^2\left(\frac{\log k}{\log K}\right)^2
- \frac{3}{4}\phi(0)^2
+\widehat{\phi}(0)\phi(0)\frac{\log k}{\log K}
-2\widehat{\phi^2}(0) \frac{\log k}{\log K}
+ O_\phi\left(\frac{\log\log k}{\log K}\right)
\nonumber \\ &\qquad
+ \textbf{P}(f;\phi)^2
+ 2\textbf{P}(f;\phi^2)
- \left(2\widehat{\phi}(0) \frac{\log k}{\log K} + \phi(0)
+ O_\phi\left(\frac{\log\log k}{\log K}\right)
\right) \textbf{P}(f;\phi).
\end{align}

From Section \ref{sec:HypothesisS}, we know that the averages of $\textbf{P}(f;\phi^2)$ and $\textbf{P}(f;\phi)$ are $o(K)$ when the support $\supp(\widehat{\phi}) \subset (-\widetilde{\sigma}/2,\widetilde{\sigma}/2)$, where $\widetilde{\sigma}$ is $\sigma(\alpha,A)$ given in Theorem \ref{thm:HypothesisS} and we assume hypothesis $H_1(\alpha,A)$. We now proceed to extracting the diagonal sum from $\textbf{P}(f;\phi)^2$ since it contributes. We have
\begin{equation}
\textbf{P}(f;\phi)^2 \ = \
\sum_{p,q} \frac{\lambda_f(p)\lambda_f(q)}{\sqrt{pq}}
\widehat{\phi}\left(\frac{\log p}{\log R}\right)
\widehat{\phi}\left(\frac{\log q}{\log R}\right)
\frac{4\log p \log q}{(\log R)^2}
\end{equation}
where $p,q$ run over the primes. When we average over $f$ we get
\begin{align}
\mathscr{P}_2(\phi) \ &\coloneqq \
\sum_{k \,{\rm even}} \frac{4\pi^2}{k-1} h\left(\frac{k-1}{K}\right)
\sum_{f\in H_k}
\frac{1}{L(1,\sym^2(f))}
\textbf{P}(f;\phi)^2
\nonumber \\ &= \
\sum_{p,q}
\widehat{\phi}\left(\frac{\log p}{\log R}\right) \widehat{\phi}\left(\frac{\log q}{\log R}\right)
\frac{4\log p \log q}{\sqrt{pq}(\log R)^2}
\mathscr{B}(p,q),
\end{align}
where $\mathscr{B}(p,q)$ was defined in (\ref{eq:mathscrB(,)}). Applying Lemma \ref{lem:ILSlemma10.1}, we have
\begin{align}
\mathscr{P}_2(\phi) \ &= \
(\widehat{h}(0)K+O(1))\sum_p \frac{1}{p} \widehat{\phi}\left(\frac{\log p}{\log R}\right)^2 \left(\frac{2 \log p}{\log R}\right)^2 + O_{\phi,h}(P^2K^{-4})
\nonumber \\ &\quad
- \pi^{1/2} K \sum_{p,q}
\frac{1}{p^{3/4} q^{3/4}}
\widehat{\phi}\left(\frac{\log p}{\log R}\right)
\widehat{\phi}\left(\frac{\log q}{\log R}\right)
\frac{4\log p \log q}{(\log R)^2}
\nonumber \\ &\qquad\qquad\qquad\quad \times
{\rm Im}\left(
\overline{\zeta}_8 \ \sum_{c=1}^\infty \
c^{-1/2} S(p,q;c) \ e\left(\frac{2\sqrt{pq}}{c}\right) \hslash\left(\frac{cK^2}{8\pi\sqrt{pq}}\right)
\right).
\end{align}
The diagonal sum yields
\begin{align}
\sum_{p} \frac{1}{p}\widehat{\phi}\left(\frac{\log p}{\log R}\right)^2
\left(\frac{2\log p}{\log R}\right)^2 \ &= \
\int_2^\infty \widehat{\phi}\left(\frac{\log x}{\log R}\right)^2
\frac{4\log^2 x}{x \log^2 R} \,d\pi(x)
\nonumber\\ &= \
\int_2^\infty \widehat{\phi}\left(\frac{\log x}{\log R}\right)^2 \frac{4\log^2 x}{x \log^2 R} \,\left(\frac{dx}{\log x} +
O\left(\frac{\log^2 x}{\sqrt{x}}\right)\,dx
\right)
\nonumber\\ &= \
\int_0^\infty \widehat{\phi}(u)^2 4u\,du + O\left(\frac{1}{\log R}\right) +
O\left(
\eta + \frac{\log R}{R^{\eta}}
\right).
\end{align}
Taking $\eta = A\log\log R / \log R$ with $A\geq 2$ a constant yields an error of size $O(\log\log R/\log R)$, where we have used the Riemann hypothesis for $\zeta(s)$ for a good error term in the integration.

We now analyze the non-diagonal component of $\mathscr{P}_2(\phi)$, denoted $\mathscr{P}_2^{\rm (nd)}(\phi)$, similar to our treatment of $\mathscr{P}(\phi)$ beginning at (\ref{eq:P(phi) bound before IBP}). We take $R \asymp K^2$.
\begin{align}
& \mathscr{P}_2^{\rm (nd)}(\phi)\nonumber\\ & \ll \
P^2K^{-4} + K\sum_{c=1}^\infty c^{-1/2} \left|
\sum_{p,q} S(p,q;c) e\left(\frac{2\sqrt{pq}}{c}\right) \hslash\left(\frac{cK^2}{8\pi\sqrt{pq}}\right)
\widehat{\phi}\left(\frac{\log p}{2\log K}\right)
\widehat{\phi}\left(\frac{\log q}{2\log K}\right)
\frac{\log p \log q}{p^{3/4}q^{3/4}(\log K)^2}
\right|
\nonumber \\ &\ll \
P^2K^{-4} + K \sum_{c=1}^\infty c^{-1/2}
\nonumber \\ &\qquad\qquad \times
\sum_{a,b \ ({\rm mod} \ c)}\hspace{-15pt}^\star
\
|S(a,b;c)|
\left| \
\underset{q \ \equiv \ b(c)}{\sum_{p \equiv a(c)}} e\left(\frac{2\sqrt{pq}}{c}\right) \hslash\left(\frac{cK^2}{8\pi\sqrt{pq}}\right)
\widehat{\phi}\left(\frac{\log p}{2\log K}\right)
\widehat{\phi}\left(\frac{\log q}{2\log K}\right)
\frac{\log p \log q}{p^{3/4}q^{3/4}(\log K)^2}
\right|.
\end{align}
We now perform summation by parts on our double sum over primes in residue classes. We let
\begin{equation}
\psi_c(x,y) \ \coloneqq \ \hslash\left(\frac{cK^2}{8\pi\sqrt{xy}}\right) \widehat{\phi}\left(\frac{\log x}{2\log K}\right) \widehat{\phi}\left(\frac{\log y}{2\log K}\right) \frac{\log x \log y}{x^{3/4}y^{3/4}(\log K)^2}
\end{equation}
and sum by parts. We have
\begin{equation}
\underset{q \ \equiv \ b(c)}{\sum_{p \equiv a(c)}} e\left(\frac{2\sqrt{pq}}{c}\right)
\psi_c(pq)
\ = \
\int_2^P \int_2^P
\mathcal{E}_2(x,y)
\frac{\partial^2 \psi_c}{\partial x \partial y}(x,y) \, dx\,dy,
\end{equation}
where
\begin{equation}
\mathcal{E}_2(x,y)
\ \coloneqq \
\underset{q \ \equiv \ b(c)}{\sum_{p \ \equiv \ a(c)}} e\left(\frac{2\sqrt{pq}}{c}\right).
\end{equation}

\emph{We now assume hypothesis $H_2(\alpha,A)$, found in (\ref{eq:HypothesisT}).}

As was illustrated by estimates in Section \ref{sec:HypothesisS}, it suffices to only integrate over the region $xy \geq c^2K^{4-2\varepsilon}$; we may also truncate our sum at $c \leq C = P K^{\varepsilon-2}$. Allowing this reduction along with Weil's estimate and $\partial^2_{xy} \psi_c(x,y) = O(K^{\varepsilon}(xy)^{-7/4})$, we obtain
\begin{align}
\mathscr{P}_2^{\rm (nd)}(\phi) \ &\ll_{\phi,h,\varepsilon} \
P^2K^{-4} + K^{1+\varepsilon} \sum_{c \ \leq \ C} c^2 \iint_{xy \ \geq \ c^2K^{4-2\varepsilon}} c^A (xy)^{\alpha-7/4+\varepsilon} \, dxdy
\nonumber \\ &\ll \
P^2K^{-4} + P^{3/2+A+2\alpha} K^{-6-2A+O(\varepsilon)}.
\end{align}
This gives an increased support of $\supp(\widehat{\phi}) \subset (-\sigma/2,\sigma/2)$ where $\sigma = 2 + \frac{6-8\alpha}{3+2A+4\alpha}$.

This concludes our proof of Theorem \ref{thm:HypothesisT}, since equation (\ref{eq:2level explicit expansion}) for the $2$-level expansion yields the main term $\langle\Phi,W_{2,\Orth}\rangle$ plus an error term that is controlled by the family averages of $\textbf{P}(f;\phi^2)$ and $\mathscr{P}_2^{\rm (nd)}(\phi)$.
Indeed, our main term is
\begin{equation}
2\int_{-\infty}^\infty |u|\widehat{\phi}(u)^2\,du + \widehat{\phi}(0)^2 - \frac{3}{4}\phi(0)^2 + \widehat{\phi}(0)\phi(0) - 2\widehat{\phi^2}(0),
\end{equation}
and it is a standard calculation (see \cite{HM, C--}) that this is equal to $\langle\Phi,W_{2,\Orth}\rangle$ (it is straightforward combinatorics to pass between $n$-level densities and $n$\textsuperscript{th} centered moments).

By Hypotheses $H_1(\alpha,A)$ and $H_2(\alpha,A)$ the error does not contribute as $K\to\infty$ for
\begin{equation}
2\sigma \ \leq \ \min\left\{\frac{5}{2}, \ 2+\frac{6-8\alpha}{1+2A+4\alpha}, \ 2+\frac{6-8\alpha}{3+2A+4\alpha}\right\} \ =\  2+\frac{6-8\alpha}{3+2A+4\alpha}.
\end{equation}
This is because the rest of the calculation for the two-level support consists of cross terms that already appear in the one-level case, and these terms are negligible assuming $H_1(\alpha,A)$ and $2\sigma$ in the above range by following the proof of Theorem \ref{thm:HypothesisS} in Section \ref{sec:HypothesisS}.


\section{An $n$-Level Conjecture}
\label{sec:n-level conjecture}

Based on our calculations in Sections \ref{sec:HypothesisS} and \ref{sec:HypothesisT}, we notice that the fundamental sums arising from local factors which may contribute have a predictable analysis. The natural generalization of cancellation in such sums lead to the following conjecture on extending support.

\begin{conj}[$n$-level Hypothesis]
\label{conj:n-level hypothesis}
Assume GRH. Let $\Phi:\R^n\to\R$ be a test function given by $\Phi(x) = \prod_{i=1}^n \phi(x_i)$.
Assume the following square-root cancellation hypothesis holds:
\begin{equation}
H_n(\alpha,A): \qquad
\underset{p_1 \ \equiv \ a_1 \ (c)}{\sum_{p_1 \ \leq \ x_1}} \ \dotsi \
\underset{p_n \ \equiv \ a_n \ (c)}{\sum_{p_2 \ \leq \ x_n}}
e\left(\frac{2\sqrt{p_1\dotsi p_n}}{c}\right)
\ \ll_\varepsilon \
c^A (x_1\dotsi x_n)^{\alpha+\varepsilon},
\end{equation}
where the implied constant is uniform over $c$ and the residue classes $a_1,\ldots,a_n$.
Then $\Phi$ agrees with random matrix theory whenever ${\rm supp}(\widehat{\phi}) \subset (-\sigma,\sigma)$, and $\sigma$ is given by
\begin{equation}
\label{eq:n-level conjecture}
n\sigma \ = \
\min\left\{
\frac{5}{2}, \
2 + \frac{6-8\alpha}{1+2A+4\alpha}, \
\ldots, \
2 + \frac{6-8\alpha}{2n-1+2A+4\alpha}
\right\} \ = \ 2 + \frac{6-8\alpha}{2n-1+2A+4\alpha}.
\end{equation}
\end{conj}

Expansions of the $n$\textsuperscript{th} centered moments or $n$-level densities have been calculated for related families where we split by sign and let the level $N\to\infty$ through the primes; see \cite{C--} for families of such cuspidal newforms, with the resulting agreement holding for ${\rm \supp}(\widehat{\phi}) \subset (-2/n, 2/n)$. A good future project\footnote{It is likely that this will be investigated in the 2023 SMALL REU.} would be to apply their combinatorial arguments in the level 1 setting, in particular seeing whether or not these are the only terms that restricted the support to $2/n$.





\appendix
\section{Theorem \ref{thm:HypothesisS} in the Literature}
\label{sec:HistoryOfThm}
We provide some brief remarks on the support improvements from Hypothesis S.

Let us suppose that $\alpha = 1/2$ and $A = 0$ for hypothesis $H_1(\alpha,A)$, i.e., the exponential sums satisfy true square-cancellation. Theorem \ref{thm:HypothesisS} says that we may take $\sigma = 5/2$. In their 1999 pre-print, Iwaniec, Luo and Sarnak perform the correct integration by parts technique, but their estimate of a sum of Bessel functions is non-sharp; they write
\begin{equation}
\label{eq:BadBesselEstimate}
\sum_{k \ {\rm even}} 2h\left(\frac{k-1}{K}\right) i^{k} J_{k-1}(x) \ = \
{\rm Main\,Term} + O_h\left(\frac{x}{K^4} + \frac{x^2}{K^6}\right),
\end{equation}
see Proposition 1 of Section 2 in their paper.
It is this extra error of $O_h(x/K^6)$ that affects the final result, giving $\sigma = 7/3$ instead of $\sigma = 5/2$. In their 2000 publication, the error term in (\ref{eq:BadBesselEstimate}) is corrected to $O_h(x/K^4)$, which leads to the correct estimate in Corollary 8.2 of \cite{ILS}. However, \cite{ILS} derives $\sigma = 22/9$ instead of $\sigma = 5/2$; this is because although their Bessel sum estimate is correct, they derive
\begin{equation}
\mathscr{P}(\phi) \ \ll_{\phi,h,\varepsilon} \
PK^{-4} + P^{\alpha+A/2+5/8}K^{-2A-9/2+\varepsilon}.
\end{equation}
It is unclear how they did this; following their integration by parts technique in the 1999 preprint, it seems that with the correction to (\ref{eq:BadBesselEstimate}) found in their 2000 publication, they would have derived equation (\ref{eq:Final Estimate for P(phi)}) as we present it today.


\section{Explicit Formulae For The Densities}
\label{sec:Explicit Formula}

We order the zeros of $L(s,f)$ by $0 \leq \gamma^{(1)} \leq \gamma^{(2)} \leq \cdots$ including multiplicity, and $\gamma^{(-j)} = -\gamma^{(j)}$, giving us a zero $\gamma^{(j)}$ for each $j\in\Z-\{0\}$.
We recall the one-level density of $L(s,f)$ with respect to a test function $\phi$ to be
\begin{equation}
D_1(f;\phi)
\ \coloneqq \
\sum_{j \neq 0} \
\phi\left(\frac{\log c_f}{2\pi}\gamma^{(j)}\right),
\end{equation}
where $c_f = k^2N$ is the analytic conductor of a form $f$ with weight $k$ and level $N$. Throughout this paper, we exclusively deal with level $1$ modular forms, so $\log c_f = 2 \log k$ for $f\in H_k$.

Similarly, the two-level density with respect to a one-dimensional test function $\phi:\R\to\R$ (so $\Phi(x_1, x_2) = \phi(x_1) \phi(x_2)$) is
\begin{equation}
D_2(f;\phi)
\ \coloneqq \
\underset{j \neq \pm \ell}{\sum_{j,\ell \neq 0}}
\phi\left(
\frac{\log c_f}{2\pi}
\gamma^{(j)}\right)
\phi\left(
\frac{\log c_f}{2\pi}
\gamma^{(\ell)}
\right).
\end{equation}
Straightforward inclusion-exclusion yields that this two-level density simplifies as
\begin{equation}
\label{eq:2level as 1levels}
D_2(f;\phi) \ = \
D_1(f;\phi)^2 - 2D_1(f;\phi^2),
\end{equation} with similar expansions for the $n$-level densities as a degree $n$ polynomial expression in terms of the one-level densities of powers of $\phi$.

In Section 4 of \cite{ILS}, they derive the following expansion. 

\begin{lem}[Lemma 4.1 of \cite{ILS}]
Let $\phi$ be even Schwartz and $\widehat{\phi}$ compactly supported. Then for $f \in H_k^{\star}(N)$, we have
\begin{align}
D_1(f;\phi) \ &= \
\widehat{\phi}(0)\frac{\log k^2N}{\log R} + \frac{1}{2}\phi(0)
- \sum_p \lambda_f(p) \widehat{\phi}\left(\frac{\log p}{\log R}\right) \frac{2 \log p}{\sqrt{p} \log R}
\nonumber \\ &\qquad
- \sum_p \lambda_f(p^2) \widehat{\phi}\left(\frac{2 \log p}{\log R}\right) \frac{2 \log p}{p \log R} + O\left(\frac{\log\log 3N}{\log R}\right)
\end{align}
where $R>1$ is an arbitrary parameter and the implied constant depends only on $\phi$.
\end{lem}

\cite{ILS} remark that, by assuming the Riemann hypothesis for $L(s,\sym^2(f))$, the second sum over primes above is a small error of size $O\left(\frac{\log\log kN}{\log R}\right)$ which depends on $\phi$, and very crucially does not depend on $f$ (cf.\ equation (4.23) of \cite{ILS}).
They also remark that one can achieve a similar error size on average via unconditional means, which is described in Appendix B of \cite{ILS}. Thus, assuming GRH, we have
\begin{equation}
\label{eq:1levelexplicit}
D_1(f;\phi) \ = \
\widehat{\phi}(0)\frac{\log k^2N}{\log R} + \frac{1}{2}\phi(0)
- \sum_p \lambda_f(p) \widehat{\phi}\left(\frac{\log p}{\log R}\right) \frac{2 \log p}{\sqrt{p} \log R}
+ O_\phi\left(\frac{\log\log kN}{\log R}\right).
\end{equation}
Let $\textbf{P}(f;\phi)$ denote the sum over primes which may contribute to the density.
\begin{equation}
\label{eq:P(f;phi)}
\textbf{P}(f;\phi)
\ \coloneqq \
\sum_p \frac{\lambda_f(p)}{\sqrt{p}}
\widehat{\phi}\left(\frac{\log p}{\log R}\right) \frac{2\log p}{\log R}.
\end{equation}
We drop the parameter $R$ from the notation of $\textbf{P}(f;\phi)$, since we often take $R$ as the analytic conductor $c_f = k^2N$ or the average conductor over the family, i.e., $R \asymp K^2N$.

We now consider the two-level density for $\phi$. By applying (\ref{eq:1levelexplicit}) to (\ref{eq:2level as 1levels}) with $f\in H_k = H_k^{\star}(1)$ and $R \asymp K^2$, we obtain
\begin{align}
\label{eq:2level explicit expansion}
D_2(f;\phi) \ &= \
\widehat{\phi}(0)^2\left(\frac{\log k}{\log K}\right)^2
- \frac{3}{4}\phi(0)^2
+\widehat{\phi}(0)\phi(0)\frac{\log k}{\log K}
-2\widehat{\phi^2}(0) \frac{\log k}{\log K}
+ O_\phi\left(\frac{\log\log k}{\log K}\right)
\nonumber \\ &\qquad
+ \textbf{P}(f;\phi)^2
+ 2\textbf{P}(f;\phi^2)
- \left(2\widehat{\phi}(0) \frac{\log k}{\log K} + \phi(0)
+ O_\phi\left(\frac{\log\log k}{\log K}\right)
\right) \textbf{P}(f;\phi).
\end{align}
As we observe in Section \ref{sec:HypothesisT}, the only prime sum which contributes to random matrix theory agreement is $\textbf{P}(f;\phi)^2$: the terms involving $\textbf{P}(f;\phi^2)$ and $\textbf{P}(f;\phi)$ become negligible under the support conditions for $\Phi(x) = \phi^2(x_1)\phi^2(x_2)$ given by Theorem \ref{thm:HypothesisS} assuming Hypothesis \ref{hyp:HypothesisS}.


\ \\

\end{document}